\documentclass[11pt]{amsart}
\usepackage{amsmath,amssymb}
\usepackage{amscd,eucal,amsthm}
\newtheorem{Theorem}{Theorem}
\newtheorem{Corollary}{Corollary}
\newtheorem{Lemma}{Lemma}
\newtheorem{Fact}{Fact}

\theoremstyle{remark}
\newtheorem*{Remark}{Remark}

\renewcommand{\to}[1][]{\xrightarrow{#1}}

\renewcommand{\gg}{{\mathfrak{g}}}
\renewcommand{\ll}{{\mathfrak{l}}}

\newcommand{\z}{{\mathfrak{z}}}
\newcommand{\s}{{\mathfrak{s}}}

\newcommand{\hh}{{\mathfrak{h}}}

\newcommand{\cc}{\mathbb{C}}

\newcommand{\A}{\mathcal{A}}
 
\DeclareMathOperator{\gr}{gr}

\DeclareMathOperator{\End}{End}

 \DeclareMathOperator{\rk}{rk}
\DeclareMathOperator{\Tr}{Tr}
\renewcommand{\phi}{\varphi}

 \makeatletter
\def\@mult#1{\raise #1\rlap{$\cdot$}\lower #1\rlap{$\cdot$}\cdot}
\def\did{\mathrel{\@mult{3pt}}}
\makeatother
\def\openrow#1#2#3{\setbox0=\vbox{\hbox
    {\vrule height#2 width#3\kern#2\vrule height#2 width0pt}\hrule height#3}
    \hbox{\leaders\copy0\hskip#1\wd0\vrule width#3}}
\def\row#1#2#3{\vbox{\hrule height#3\openrow{#1}{#2}{#3}}}
\def\Yr#1{\row{#1}{1.5ex}{.1ex}}

\def\DY#1\endDY{\baselineskip=1ex\lineskip=0pt\lineskiplimit=0pt{\vcenter
    {\Yr#1}}}
\def\openclm#1#2#3{\setbox0=\vbox{\hrule height#3\hbox
    {\vrule width0pt\kern#2\vrule width#3 height#2}}\vtop
    {\leaders\copy0\vskip#1\ht0\hrule height#3}}
\def\clm#1#2#3{\hbox{\vrule width#3\openclm{#1}{#2}{#3}}}
\def\Yc#1{\clm{#1}{1.5ex}{.1ex}}

\def\CDY#1\endCDY{{\vcenter{\hbox{\Yc#1}}}}

\author{L.~G.~Rybnikov}
\title{Uniqueness of higher Gaudin hamiltonians}
\address{Poncelet laboratory (Independent University of Moscow and CNRS) and Moscow State University,
department of Mechanics and Mathematics}
\email{leo.rybnikov@gmail.com}
\thanks{The work was partially supported by RFBR grant 04-01-00702, RFBR grant 05
01 00988-a and RFBR grant 05-01-02805-CNRSL-a.}

\begin{document}
\maketitle
\begin{abstract}
For any semisimple Lie algebra $\mathfrak{g}$, the universal
enveloping algebra of the infinite-dimensional pro-nilpotent Lie
algebra $\mathfrak{g}_-:=\mathfrak{g}\otimes t^{-1}\cc[t^{-1}]$
contains a large commutative subalgebra $\mathcal{A}\subset
U(\mathfrak{g}_-)$. This subalgebra comes from the center of the
universal enveloping of the affine Kac--Moody algebra
$\hat{\mathfrak{g}}$ at the critical level and gives rise to the
construction of higher hamiltonians of the Gaudin model (due to
Feigin, Frenkel and Reshetikhin). Though there are no explicit
formulas for the generators of $\mathcal{A}$ known in general, the
"classical analogue" of this subalgebra, i.e. the associated
graded subalgebra in the Poisson algebra $A\subset
S(\mathfrak{g}_-)$, can be easily described. In this note we show
that the "classical" subalgebra $A\subset S(\mathfrak{g}_-)$ is
the Poisson centralizer of some of its quadratic elements, and
deduce from this that the "quantum" subalgebra $\mathcal{A}\subset
U(\mathfrak{g}_-)$ is uniquely determined by the space of
quadratic elements of the classical one. In particular, this means
that some different constructions of higher Gaudin hamiltonians
(namely, Feigin-Frenkel-Reshetikhin's method and
Talalaev-Chervov's method), give the same family of commuting
operators. The proof uses some ideas of the previous paper
math.QA/0608586.
\end{abstract}

\section{Introduction.} Let $\gg$ be a semisimple complex Lie
algebra. We consider the infinite-dimensional pro-nilpotent Lie
algebra $\gg_-:=\gg\otimes t^{-1}\cc[t^{-1}]$ -- it is a "half" of
the corresponding affine Kac--Moody algebra $\hat\gg$. The
universal enveloping algebra $U(\gg_-)$ bears a natural filtration
by the degree with respect to the generators. The associated
graded algebra is the symmetric algebra $S(\gg_-)$ by the
Poincar\'e--Birkhoff--Witt theorem. The commutator operation on
$U(\gg_-)$  defines the Poisson--Lie bracket $\{\cdot,\cdot\}$ on
$S(\gg_-)$: for the generators $x,y\in\gg_-$ we have
$\{x,y\}=[x,y]$. For any $g\in\gg$, we denote the element
$g\otimes t^{-m}\in\gg_-$ by $g[-m]$.

The Poisson algebra $S(\gg_-)$ contains a large
Poisson-commutative subalgebra $A\subset S(\gg_-)$. This
subalgebra can be constructed as follows.

Consider the following derivation of the Lie algebra $\gg_-$:
\begin{equation}\label{der1}
\partial_t(g[-m])=-mg[-m-1]\quad\forall g\in\gg, m=-1,-2,\dots
\end{equation}
The derivation (\ref{der1}) extends to the derivation of the
associative algebras $S(\gg_-)$ and $U(\gg_-)$.

Let $i_{-1}:S(\gg)\hookrightarrow S(\gg_-)$ be the embedding,
which maps $g\in\gg$ to $g[-1]$. Let $\Phi_k,\ k=1,\dots,\rk\gg$
be the generators of the algebra of invariants $S(\gg)^{\gg}$.

\begin{Fact}\cite{BD,FFR,Fr2} \begin{enumerate} \item The subalgebra  $A\subset S(\gg_-)$ generated by the elements
$\partial_t^n \overline{S_k}$, $k=1,\dots,\rk\gg$,
$n=0,1,2,\dots$, where $\overline{S_k}=i_{-1}(\Phi_k)$, is
Poisson-commutative. \item There exist the elements $S_k\in
U(\gg_-)^\gg$ such that $\gr S_k=\overline{S_k}$ and the
subalgebra $\A\subset U(\gg_-)$ generated by $\partial_t^n S_k$,
$k=1,\dots,\rk\gg$, $n=0,1,2,\dots$ is commutative.\end{enumerate}
\end{Fact}

\begin{Remark} The generators of the subalgebra $A\subset S(\gg_-)$ can be described in the following
equivalent way. Let $i(z):S(\gg)\hookrightarrow S(\gg_-)$ be be
the embedding depending on the formal parameter $z$, which maps
$g\in\gg$ to $\sum\limits_{k=1}^{\infty}z^{k-1}g[-k]$. Then the
coefficients of the power series $\overline{S_k}(z)=i(z)(\Phi_k)$
in $z$ freely generate the subalgebra $A\subset S(\gg_-)$.
\end{Remark}

\begin{Remark} The subalgebra $\A\subset U(\gg_-)$ comes from the center of
$U(\hat\gg)$ at the critical level by the AKS-scheme (see
\cite{FFR,ER,ChT}).
\end{Remark}

\begin{Remark} No general explicit formulas for the elements
$S_k$ are known in general. For the quadratic Casimir element
$\Phi_1$, the corresponding element $S_1\in\A$ is obtained from
$\overline{S_1}=i_{-1}(\Phi_1)$ by the symmetrization map. More
explicitly, we have $S_1=\sum\limits_{a=1}^{\dim\gg}x_a[-1]^2$,
where $\{x_a\},\ a=1,\dots,\dim\gg$ is an orthonormal (with
respect to the Killing form) basis of $\gg$.
\end{Remark}

\begin{Remark} For
$\gg=sl_r$ explicit formulas for $S_k$ were obtained by Talalaev
and Chervov in \cite{ChT}. We shall reproduce these formulas in
section~2.
\end{Remark}

The main result of the present paper is the following

\begin{Theorem}\label{main} \begin{enumerate}
\item The subalgebra $A\subset S(\gg_-)$ is a Poisson centralizer
of the element $\overline{S_1}=i_{-1}(\Phi_1)\in A$ (where
$\Phi_1$ is the quadratic Casimir). In particular, $A\subset
S(\gg_-)$ is a \emph{maximal} Poisson-commutative subalgebra in
$S(\gg_-)$. \item The subalgebra $\A\subset U(\gg_-)$ is a
centralizer of the element $S_1\in\A$. \item Any commutative
subalgebra $\tilde\A\subset U(\gg_-)^\gg$, such that
$\gr\tilde\A=A$, coincides with $\A\subset U(\gg_-)^\gg$ (in other
words, there is a unique $G$-invariant lifting of the commutative
subalgebra $A\subset S(\gg_-)$ to the universal enveloping
algebra).
\end{enumerate}
\end{Theorem}

\begin{Remark} In \cite{R} a similar fact for Mischenko--Fomenko
subalgebras of the Poisson algebra $S(\gg)$ is proved.
Unfortunately, the ideas of \cite{R} can not be generalized
straightforwardly to our situation. The proof of
Theorem~\ref{main} uses also some new ideas.
\end{Remark}

We will prove the Theorem in section~3. Theorem~\ref{main} has
some applications in the quantization of higher hamiltonians of
the Gaudin model. We discuss this in section~2.

I thank B.~L.~Feigin, A.~V.~Chervov, and D.~V.~Talalaev for useful
discussions.

\section{Gaudin model} Gaudin model was introduced in \cite{G1} as a spin
model related to the Lie algebra $sl_2$, and generalized to the
case of an arbitrary semisimple Lie algebra in \cite{G}, 13.2.2.
The generalized Gaudin model has the following algebraic
interpretation. For any $x\in\gg$, set
$x^{(i)}=1\otimes\dots\otimes 1\otimes x\otimes
1\otimes\dots\otimes 1\in U(\gg)^{\otimes n}$ ($x$ stands at the
$i$th place). Let $\{x_a\},\ a=1,\dots,\dim\gg$, be an orthonormal
basis of $\gg$ with respect to Killing form, and let
$z_1,\dots,z_n$ be pairwise distinct complex numbers. The
hamiltonians of Gaudin model are the following commuting elements
of $U(\gg)^{\otimes n}$:
\begin{equation}\label{quadratic}
H_i=\sum\limits_{k\neq i}\sum\limits_{a=1}^{\dim\gg}
\frac{x_a^{(i)}x_a^{(k)}}{z_i-z_k}.
\end{equation}

In \cite{FFR}, a large commutative subalgebra
$\A(z_1,\dots,z_n)\subset U(\gg)^{\otimes n}$ containing $H_i$ was
constructed. Generators of this algebra are known as higher Gaudin
hamiltonians. The commutative subalgebra $\A(z_1,\dots,z_n)\subset
U(\gg)^{\otimes n}$ is the image of the subalgebra $\A\subset
U(\gg_-)$ under the homomorphism $U(\gg_-)\to U(\gg)^{\otimes n}$
of specialization at the points $z_1,\dots,z_n$ (see
\cite{FFR,ER}).

In~\cite{Tal} D.~Talalaev constructed explicitly some elements of
$U(\gg)^{\otimes n}$ commuting with quadratic Gaudin hamiltonians
for the case $\gg=\gg\ll_r$. The formulas of~\cite{Tal} are
universal, i.e. actually they describe a commutative subalgebra of
$U(\gg_-)$ which gives a commutative subalgebra of
$U(\gg)^{\otimes n}$ as the image of the specialization
homomorphism at the points $z_1,\dots,z_n$ (see \cite{ChT}).

Namely, set
$$L(z)=\sum\limits_{1\le i,j\le
r}\sum\limits_{n=1}^{\infty}z^{n-1}e_{ij}[-n]\otimes e_{ji}\in
U(\gg_-)\otimes\End\cc^r,$$ where $z$ is a formal parameter, and
consider the following differential operator in $z$ with the
coefficients from $U(\gg_-)$:
$$
D=\Tr A_r\prod\limits_{i=1}^r(L(z)^{(i)}-\partial_z)=
\partial_z^r+\sum\limits_{k=1}^{r}\sum\limits_{n=1}^{\infty}Q_{n,k}z^{k-1}\partial_z^{r-k}.
$$

Here we denote by $L(z)^{(i)}\in
U(\gg_-)\otimes(\End\cc^r)^{\otimes r}$ the element obtained by
putting $L(z)$ to the $i$-th tensor factor, and $A_r$ denotes the
projector onto $U(\gg_-)\otimes\End(\Lambda^r\cc^r)\subset
U(\gg_-)\otimes(\End\cc^r)^{\otimes r}$. It follows from the works
\cite{ChT,Tal} by Chervov and Talalaev that the elements
$Q_{n,k}\in U(\gg_-)$ pairwise commute.

It is easy to see the Poisson-commutative subalgebra in $S(\gg_-)$
generated by $\gr Q_{n,k}\in S(\gg_-)$ coincides with $A$, i.e.
$\gr
\sum\limits_{n=1}^{\infty}Q_{n,k}z^{k-1}=\overline{S_k}(z)=i(z)(\Phi_k)\in
S(\gg_-)$, where we take as the generators $\Phi_k$ of the algebra
of invariants $S(\gg)^\gg$ the coefficients of the characteristic
polynomial (see Remark~1 in \cite{ChT}).

From Theorem~\ref{main} we obtain the following

\begin{Corollary}\label{Dima} In the case of $\gg=\gg\ll_r$, the higher Gaudin hamiltonians in the sense
of~\cite{Tal} coincide with higher Gaudin hamiltonians in the
sense of~\cite{FFR}.
\end{Corollary}

\section{Proof of Theorem~\ref{main}.}
We fix a Cartan subalgebra $\hh\subset\gg$. Let $e,\ h,\ f$ be a
principal $\s\ll_2$-triple in the Lie algebra $\gg$ (here
$h\in\hh$). Let $\z_{\gg}(f)$ be the centralizer of $f$ in $\gg$.
Now, we will prove some analogs of the following classical results
for the subalgebra $A\subset S(\gg_-)$.

\begin{Fact}\label{kostant}\emph{(Kostant \cite{Ko})} The homomorphism of the restriction to
the affine subspace $\cc[\gg]\to \cc[e+\z_{\gg}(f)]$ induces an
algebra isomorphism $\cc[\gg]^{\gg}\tilde\to \cc[e+\z_{\gg}(f)]$.
\end{Fact}

\begin{Fact}\label{chevalle}\emph{(Chevalley, see e.g. \cite{He})} The image of $\cc[\gg]^{\gg}$ under the restriction to
the Cartan subalgebra $\hh\subset\gg$ is the algebra of invariants
of the Weyl group action on $\hh$. In particular $\cc[\hh]$ is an
algebraic extension of the image of $\cc[\gg]^{\gg}$.
\end{Fact}

Let $V\subset\gg$ be an $h$-invariant complement of the subspace
$\z_{\gg}(f)\subset\gg$. We consider the homomorphism
$\pi:S(\gg_-)\to S(\z_{\gg}(f)_-)$ defined as follows.

$$
\pi(x[-m])=x[-m]\quad\forall x\in\z_{\gg}(f),
$$
and
$$
\pi(x[-m])=\delta_{1m}\langle x,f\rangle\quad\forall x\in V.
$$

Let $\psi:\gg_-\to\hh_-$ be an $\hh$-invariant projection. This
projection extends to the algebra homomorphism $S(\gg_-)\to
S(\hh_-)$, which we denote by the same letter $\psi$.

\begin{Lemma}\label{sl-triple}\begin{enumerate} \item
The homomorphism $\pi:A\to S(\z_{\gg}(f)_-)$ is an isomorphism of
algebras. \item $\psi(A)\subset S(\hh_-)$ is an algebraic
extension.
\end{enumerate}
\end{Lemma}
\begin{proof} Let us prove the first assertion.
By Fact~\ref{kostant} the restriction of the homomorphism $\pi$ to
the subalgebra $\cc[\overline{S_1},\dots,\overline{S_l}]\subset A$
induces an algebra isomorphism
$\pi:\cc[\overline{S_1},\dots,\overline{S_l}]\to
S(\z_{\gg}(f)[-1])$. Now, it remains to notice that the
homomorphism $\pi:S(\gg_-)\to S(\z_{\gg}(f)_-)$ commutes with the
operator $\partial_t$ and therefore
$$\pi:A=\cc[\overline{S_1},\dots,\overline{S_l},\dots,\partial_t^k\overline{S_1},\dots,
\partial_t^k\overline{S_l},\dots]\to S(\z_{\gg}(f)_-)$$ is an algebra isomorphism as well.

Now, let us prove the second assertion of the Lemma. By
Fact~\ref{chevalle} we have that every element of the subspace
$\hh[-1]\subset S(\hh_-)$ is algebraic over $\psi(A)$. Since the
homomorphism $\psi:S(\gg_-)\to S(\hh_-)$ commutes with the
operator $\partial_t$ and the subalgebra $A\subset S(\gg_-)$ is
stable under $\partial_t$, we see that for any positive integer
$k$ every element of the subspace $\partial_t^k(\hh[-1])=\hh[-k]$
is algebraic over $\psi(A)$ as well. Thus, $\psi(A)\subset
S(\hh_-)=S(\bigoplus\limits_{k=1}^{\infty}\hh[-k])$ is an
algebraic extension.
\end{proof}

\begin{Lemma}\label{alg_closed} The subalgebra
$A$ is algebraically closed in $S(\gg_-)$.
\end{Lemma}
\begin{proof}
Assume that there is an element $a\in S(\gg_-)$ such that
$a\not\in A$ ¨ $p_Na^N+p_{N-1}a^{N-1}+\dots+p_0=0$ for some
$p_0,\dots,p_N\in A$. Suppose that the number $N$ is minimal
possible. By the first assertion of Lemma~\ref{sl-triple} without
loss of generality we can assume that $\pi(a)=0$ without loss of
generality. This means, in particular, that $\pi(p_0)=0$, and
hence, by the first assertion of Lemma~\ref{sl-triple} $p_0=0$.
Since $S(\gg_-)$ has no divisors of zero, we have
$p_Na^{N-1}+p_{N-1}a^{N-2}+\dots+p_1=0$. This means that $N$ is
not minimal.
\end{proof}

Now, we define a $1$-parametric family of automorphisms $\phi_s$
of the Poisson algebra $S(\gg_-)$: for $x\in\gg$ set
$\phi_s(x[-k])=x[-k]+s\delta_{1,k}\langle h,x\rangle$ (here
$\langle \cdot,\cdot\rangle$ denotes the Killing form). The
following assertion is checked by direct computation.

\begin{Lemma}\label{predel}
$\lim\limits_{s\to\infty}\frac{\overline{S_1}-s^2\Phi_1(h)}{2s}=h[-1]$.
\end{Lemma}

\begin{Lemma}\label{ht}
 $S(\hh_-)$ is the centralizer of $h[-1]$ in $S(\gg_-)$.
\end{Lemma}
\begin{proof}
Fix any ordering on the set of roots of $\gg$ and define an
ordering on the generators of $S(\gg_-)$ in the following way:
$$h_{\alpha_i}[-k]<e_{\alpha}[-m]\,\forall
k,m,\alpha_i,\alpha;$$
$$e_{\alpha}[-k]<e_{\beta}[-m]\,\forall
k,m,\,\text{if}\, \alpha<\beta;
$$
$$h_{\alpha_i}[-k]<h_{\alpha_j}[-m]\,\forall
k,m,\,\text{if}\, \alpha_i<\alpha_j;
$$
$$x[-k]<x[-m],\,\text{if}\, k<m;
$$

Assume that there exists $f\in S(\gg_-)$ commuting with $h[-1]$
such that $f\not\in S(\hh_-)$. Let
$$\prod\limits_{i=1}^{l}\prod\limits_{k=1}^{\infty}(h_{\alpha_i}[-k])^{m_{i,k}}\prod\limits_{\alpha\in\Delta}
\prod\limits_{k=1}^{\infty} (e_{\alpha}[-k])^{n_{\alpha,k}}$$ be
the leading monomial of $f$ with respect to our ordering, and let
$\alpha_{max}$ be the maximal root such that $n_{\alpha,k}\ne0$
for some $k$. Let $k_{max}$ be the largest of such $k$ (by the
assumption $f\not\in S(\hh_-)$, therefore such $\alpha_{max}$ and
$k_{max}$ exist). Then the leading monomial of the commutator
$\{h[-1],f\}$ is
$$n_{\alpha_{max},k_{max}}\langle h,\alpha_{max}\rangle\frac{e_{\alpha_{max}}[-k_{max}-1]}{e_{\alpha_{max}}[-k_{max}]}
\prod\limits_{i=1}^{l}\prod\limits_{k=1}^{\infty}(h_{\alpha_i}[-k])^{m_{i,k}}\prod\limits_{\alpha\in\Delta}
\prod\limits_{k=1}^{\infty} (e_{\alpha}[-k])^{n_{\alpha,k}}\ne0.$$
This contradicts our assumption. The Lemma is proved.
\end{proof}

\emph{Proof of Theorem~\ref{main}.} (1). Let $I\subset S(\gg_-)$
be the kernel of the homomorphism $\psi: S(\gg_-)\to S(\hh_-)$,
i.e. the ideal generated by all $e_\alpha[-m],\ \alpha\in\Delta,\
m=1,2,\dots$. Note that all the automorphisms $\phi_s$ preserve
the ideal $I$. Assume that there is an element $a$ in the Poisson
centralizer of $\overline{S_1}$ in $S(\gg_-)$, which do not belong
to $A$. By Lemma~\ref{alg_closed} the element $a$ is
transcendental over $A$. By the second assertion of
Lemma~\ref{sl-triple}, we can pass to a polynomial expression in
$a$ with the coefficients from $A$ to make $\psi(a)=0$, i.e. $a\in
I$. For appropriate $k$, there exist a non-zero limit $\overline
a:=\lim\limits_{s\to\infty}\frac{\phi_s(a)}{s^k}\in I$. Moreover,
by Lemma~\ref{predel} $\overline a$ belongs to the centralizer of
$h[-1]$ in $S(\gg_-)$. On the other hand, $\overline a\in I$, an
hence $\overline a\not\in S(\hh_-)$. This contradicts
Lemma~\ref{ht}. Thus, the Poisson centralizer of $\overline{S_1}$
in $S(\gg_-)$ coincides with $A$.

(2). Let $Z$ be the centralizer of $S_1$ in the associative
algebra $U(\gg)_-$. We have $\gr Z\subset A$ and $\A\subset Z$,
hence $Z=\A$.

(3). It is clear that $S_1$ is the unique -- up to an additive
constant -- lifting of $\overline{S_1}$ to $U(\gg_-)^{\gg}$. Since
the centralizer does not depend on adding a constant, we see that
$\A$ is the unique lifting of the Poisson commutative subalgebra
$A\subset S(\gg_-)$ to the associative algebra $U(\gg_-)^{\gg}$.
Thus, Theorem~\ref{main} is proved. \qquad$\Box$


\begin{thebibliography}{10}

\bibitem[BD]{BD} Beilinson, A. and Drinfeld, V. {\em Quantization of Hitchin's
integrable system and Hecke eigen-sheaves}, Preprint, available at
www.ma.utexas.edu/$\sim$benzvi/BD.

\bibitem[ChT]{ChT} Chervov, A. and Talalaev, D.
\emph{Quantum spectral curves, quantum integrable systems and the
geometric Langlands correspondence}, preprint hep-th/0604128

\bibitem[ER]{ER} Enriquez, B. and Rubtsov, V. {\em Hitchin systems, higher Gaudin hamiltonians and $r$-matrices.}
Preprint math.alg-geom/9503010.

\bibitem[FF]{FF} Feigin, B. and Frenkel, E. {\em Affine Kac--Moody
algebras at the critical level and Gelfand--Dikii algebras}, Int.
Jour. Mod. Phys. {\bf A7}, Supplement 1A (1992) 197--215.

\bibitem[Fr1]{Fr1}
Frenkel, E. {\em Affine Algebras, Langlands Duality and Bethe
Ansatz}, q-alg/9506003

\bibitem[Fr2]{Fr2}
Frenkel,  E. {\em Lectures on Wakimoto modules, opers and the
center at the critical level}, QA/0210029.

\bibitem[FFR]{FFR} Feigin, B., Frenkel, E. and Reshetikhin, N. {\em Gaudin model, Bethe Ansatz and critical level.}
Comm. Math. Phys., 166 (1994), pp. 27-62.

\bibitem[G1]{G1}
Gaudin, M. {\em Diagonalisation d'une classe d'hamiltoniens de
spin}, J. de Physique, t.37, N 10, p. 1087--1098, 1976.

\bibitem[G]{G}Gaudin, M. \emph{La fonction d'onde de Bethe.}
 Collection du Commissariat a` l'E'nergie
Atomique: Se'rie Scientifique. Masson, Paris, 1983. xvi+331 pp.

\bibitem[He]{He} Helgason S., {\em Differential geometry and symmetric spaces,
Academic Press}, New York.

\bibitem[Ko]{Ko}Kostant, B. \emph{Lie group representations on polynomial
rings.}
Amer. J. Math.  85  1963 327--404.

\bibitem[R]{R} Rybnikov, L. G. \emph{Centralizers of certain quadratic elements in
Poisson-Lie algebras and Argument Shift method. (Russian)} Uspekhi
Mat. Nauk 60 (2005), no. 2(362), 173--174; math.QA/0608586.

\bibitem[Tal]{Tal} Talalaev, D. {\em Quantization of the Gaudin system},
Preprint hep-th/0404153.

\end{thebibliography}
\end{document}